# Harmony Search as a Metaheuristic Algorithm

**Xin-She Yang**[1],

**Abstract** This first chapter intends to review and analyze the powerful new Harmony Search (HS) algorithm in the context of metaheuristic algorithms. I will first outline the fundamental steps of Harmony Search, and how it works. I then try to identify the characteristics of metaheuristics and analyze why HS is a good metaheuristic algorithm. I then review briefly other popular metaheuristics such as particle swarm optimization so as to find their similarities and differences from HS. Finally, I will discuss the ways to improve and develop new variants of HS, and make suggestions for further research including open questions.

## 1 Introduction

When listening to a beautiful piece of classical music, who has ever wondered if there is any connection between music and finding an optimal solution to a tough design problem such as the water distribution networks or other design problems in engineering? Now for the first time ever, scientists have found such an interesting connection by developing a new algorithm, called Harmony Search. Harmony Search (HS) was first developed by Zong Woo Geem et al. in 2001 [1], though it is a relatively new metaheuristic algorithm, its effectiveness and advantages have been demonstrated in various applications. Since its first appearance in 2001, it has been applied to solve many optimization problems including function optimization, engineering optimization [2], water distribution networks [3], groundwater modelling, energy-saving dispatch, truss design, vehicle routing, and others. The possibility of combining harmony search with other algorithms such as Particle Swarm Optimization (PSO) has also been investigated.

    Harmony search is a music-based metaheuristic optimization algorithm. It was inspired by the observation that the aim of music is to search for a perfect state of harmony. This harmony in music is analogous to find the optimality in an optimization process. The search process in optimization can be compared to a

---

[1] Department of Engineering, University of Cambridge, Trumpington Street, Cambridge CB2 1PZ, UK. Email: xy227@cam.ac.uk

jazz musician's improvisation process. On the one hand, the perfectly pleasing harmony is determined by the audio aesthetic standard. A musician always intends to produce a piece of music with perfect harmony. On the other hand, an optimal solution to an optimization problem should be the best solution available to the problem under the given objectives and limited by constraints. Both processes intend to produce the best or optimum.

Such similarities between two processes can be used to develop new algorithms by learning from each other. Harmony Search is just such a successful example by transforming the qualitative improvisation process into some quantitative rules by idealization, and thus turning the beauty and harmony of music into an optimization procedure through search for a perfect harmony, namely, the Harmony Search (HS) or Harmony Search algorithm.

## 2 Harmony Search as a Metaheuristic Method

Before we introduce the fundamentals of HS algorithm, let us first briefly describe the way to describe the aesthetic quality of music. Then, we will discuss the pseudo code of HS algorithm and two simple examples to demonstrate how it works.

### *2.1 Aesthetic Quality of Music*

The aesthetic quality of a musical instrument is essentially determined by its pitch (or frequency), timbre (or sound quality), and amplitude (or loudness). Timbre is largely determined by the harmonic content that is in turn determined by the waveforms or modulations of the sound signal. However, the harmonics that it can generate will largely depend on the pitch or frequency range of the particular instrument.

Different notes have different frequencies. For example, the note A above middle C (or standard concert A4) has a fundamental frequency of $f_0=440$ Hz. As the speed of sound in dry air is about $v=331+0.6T$ m/s where $T$ is the temperature in degrees Celsius near T=0. So at room temperature T=20°C, the A4 note has a wavelength $\lambda = v/f_0 \approx 0.7795$ m. When we adjust the pitch, we are in fact trying to change the frequency. In music theory, pitch $p_n$ in MIDI is often represented as a numerical scale (a linear pitch space) using the following formula

$$p_n = 69 + 12\log_2(\frac{f}{440Hz}), \qquad (1)$$

or

$$f = 440 \times 2^{(p_n-69)/12}, \qquad (2)$$

which means that the A4 notes has a pitch number 69. On this scale, octaves correspond to size 12 while semitone corresponds to size 1, which leads to that the ratio of frequencies of two notes that are an octave apart is 2:1. Thus, the frequency of a note is doubled (halved) when it raised (lowered) an octave. For example, A2 has a frequency of 110Hz while A5 has a frequency of 880Hz.

The measurement of harmony where different pitches occur simultaneously, like any aesthetic quality, is subjective to some extent. However, it is possible to use some standard estimation for harmony. The frequency ratio, pioneered by ancient Greek mathematician Pythagoras, is a good way for such estimation. For example, the octave with a ratio of 1:2 sounds pleasant when playing together, so are the notes with a ratio of 2:3. However, it is unlikely for any random notes played by a monkey to produce a pleasant harmony.

## *2.2 Harmony Search*

In order to explain the Harmony Search in more detail, let us first idealize the improvisation process by a skilled musician. When a musician is improvising, he or she has three possible choices: (1) play any famous piece of music (a series of pitches in harmony) exactly from his or her memory; (2) play something similar to a known piece (thus adjusting the pitch slightly); or (3) compose new or random notes. Zong Woo Geem et al. formalized these three options into quantitative optimization process in 2001, and the three corresponding components become: usage of harmony memory, pitch adjusting, and randomization [1].

The usage of harmony memory is important, as it is similar to the choice of the best-fit individuals in genetic algorithms (GA). This will ensure that the best harmonies will be carried over to the new harmony memory. In order to use this memory more effectively, it is typically assigned as a parameter $r_{accept} \in [0,1]$, called harmony memory accepting or considering rate. If this rate is too low, only few best harmonies are selected and it may converge too slowly. If this rate is extremely high (near 1), almost all the harmonies are used in the harmony memory, then other harmonies are not explored well, leading to potentially wrong solutions. Therefore, typically, we use $r_{accept}=0.7 \sim 0.95$.

The second component is the pitch adjustment determined by a pitch bandwidth $b_{range}$ and a pitch adjusting rate $r_{pa}$. Though in music, pitch adjustment means to change the frequencies, it corresponds to generate a slightly different solution in the Harmony Search algorithm [1]. In theory, the pitch can be adjusted linearly or nonlinearly, but in practice, linear adjustment is used. So we have

$$x_{new} = x_{old} + b_{range} * \varepsilon \qquad (3)$$

where $x_{old}$ is the existing pitch or solution from the harmony memory, and $x_{new}$ is the new pitch after the pitch adjusting action. This essentially produces a new solution around the existing quality solution by varying the pitch slightly by a small random amount [1,2]. Here $\varepsilon$ is a random number generator in the range of [-1,1]. Pitch adjustment is similar to the mutation operator in genetic algorithms. We can assign a pitch-adjusting rate ($r_{pa}$) to control the degree of the adjustment. A low pitch adjusting rate with a narrow bandwidth can slow down the convergence of HS because the limitation in the exploration of only a small subspace of the whole search space. On the other hand, a very high pitch-adjusting rate with a wide bandwidth may cause the solution to scatter around some potential optima as in a random search. Thus, we usually use $r_{pa}$=0.1 ~0.5 in most applications.

| Harmony Search |
|---|
| **begin** |
|     *Objective function f(**x**), **x**=($x_1,x_2, …,x_d$)$^T$* |
|     *Generate initial harmonics (real number arrays)* |
|     *Define pitch adjusting rate ($r_{pa}$), pitch limits and bandwidth* |
|     *Define harmony memory accepting rate ($r_{accept}$)* |
|     **while** *( t<Max number of iterations )* |
|         *Generate new harmonics by accepting best harmonics* |
|         *Adjust pitch to get new harmonics (solutions)* |
|         **if** *(rand>$r_{accept}$), choose an existing harmonic randomly* |
|         **else if** *(rand>$r_{pa}$), adjust the pitch randomly within limits* |
|         **else**   *generate new harmonics via randomization* |
|         **end if** |
|         *Accept the new harmonics (solutions) if better* |
|     **end while** |
|     *Find the current best solutions* |
| **end** |

Figure 1: Pseudo code of the Harmony Search algorithm.

The third component is the randomization, which is to increase the diversity of the solutions. Although adjusting pitch has a similar role, but it is limited to certain local pitch adjustment and thus corresponds to a local search. The use of randomization can drive the system further to explore various diverse solutions so as to find the global optimality.

The three components in harmony search can be summarized as the pseudo code shown in Fig. 1. In this pseudo code, we can see that the probability of randomization is

$$P_{\text{random}} = 1 - r_{\text{accept}}, \tag{4}$$

and the actual probability of adjusting pitches is

$$P_{\text{pitch}} = r_{\text{accept}} * r_{\text{pa}}. \tag{5}$$

## 2.3 Implementation

The three components of the HS algorithm described in the above section can easily be implemented using any programming language, though it should be straightforward to carry out simulations with visualization using Matlab.

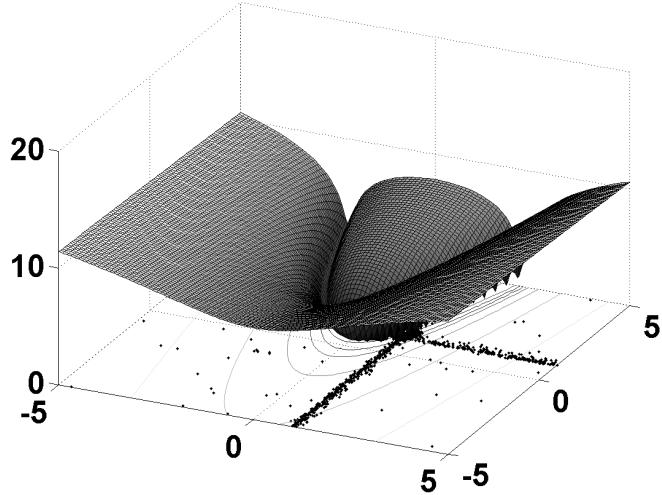

Figure 2: The search paths of finding the global optimal solution (1,1) using the harmony search.

For Rosenbrock's logarithmic banana function
$$f(x, y) = \ln[1 + (1 - x)^2 + 100(y - x^2)^2] \qquad (6)$$
where $(x,y) \in [-10,10] \times [-10,10]$, it has a global minimum $f_{min}=0$ at (1,1). The best estimate solution (1.0023,1.0070) is obtained after 15,000 iterations using the Matlab program described in [4]. In a modern PC (say, with a 3GHz processor), it usually takes about 2 minutes. The variations of these solutions and their paths are shown in Fig. 2.

We have used 20 harmonics, the harmony accepting rate $r_{accept}=0.95$, and the pitch adjusting rate $r_{pa}=0.7$. The search paths are plotted together with the landscape of $f(x,y)$. From Fig. 2, we can see that the pitch adjustment is more intensive in local regions (two thin strips), this is probably another reason why the harmony search is more efficient than genetic algorithms.

As a further example, Michalewicz's bivariate function
$$f(x, y) = -\sin(x)\sin^{20}(\frac{x^2}{\pi}) - \sin(y)\sin^{20}(\frac{2y^2}{\pi}), \qquad (7)$$

has a global minimum $f_{min} \approx -1.801$ at [2.20319,1.57049] in the domain $0 \leq x \leq \pi$ and $0 \leq y \leq \pi$. This global minimum can be found after about 23,000 function evaluations, and the results are shown in Fig. 3.

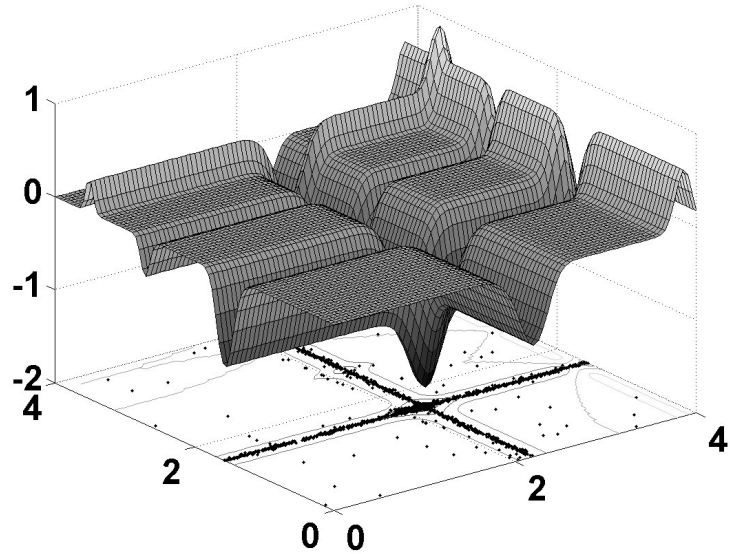

Fig.3 Harmony search for Michalewicz's bivariate function.

In addition to the above-mentioned two benchmark examples, this book will contain many successful examples of using the HS algorithm to solve various tough optimization problems, and comparing the HS results with those of other algorithms. Such a comparison for different types of algorithms and problems is still an area of active research.

## 3 Other Metaheuristics

### *3.1 Metaheuristics*

Heuristic algorithms typically intend to find a good solution to an optimization problem by 'trial-and-error' in a reasonable amount of computing time. Here 'heuristic' means to 'find' or 'search' by trials and errors. There is no guarantee to find the best or optimal solution, though it might be a better or improved solution than an educated guess. Any reasonably good solution, often suboptimal or near optimal, would be good enough for such problems. Broadly speaking, local search methods are heuristic methods because their parameter search is focused on the local variations, and the optimal or best solution can be well outside this local region. However, a high-quality feasible solution in the local region of interest is usually accepted as a good solution in many optimization problems in practice if time is the major constraint.

Metaheuristic algorithms are higher-level heuristic algorithms. Here, '*meta-*' means 'higher-level' or 'beyond', so metaheuristic means literally to find the solution using higher-level techniques, though certain trial-and-error processes are still used. Broadly speaking, metaheuristics are considered as higher-level techniques or strategies which intend to combine lower-level techniques and tactics for exploration and exploitation of the huge space for parameter search. In recent years, the word '*metaheuristics*' refers to all modern higher-level algorithms [4], including Particle Swarm Optimization (PSO), Simulated Annealing (SA), Evolutionary Algorithms (EA) including Genetic Algorithms (GA), Tabu Search (TS), Ant Colony Optimization (ACO), Bee Algorithms (BA), Firefly Algorithms (FA), and, certainly Harmony Search (HS).

There are two important components in modern metaheuristics, and they are: intensification and diversification, and such terminologies are derived from Tabu search [5]. For an algorithm to be efficient and effective, it must be able to generate a diverse range of solutions including the potentially optimal solutions so as to explore the whole search space effectively, while it intensifies its search around the neibourhood of an optimal or nearly optimal solution. In order to do so, every part of the search space must be accessible though not necessarily visited during the search. Diversification is often in the form of randomization with a random component attached to a deterministic component in order to explore the search space effectively and efficiently, while intensification is the exploitation of past solutions so as to select the potentially good solutions via elitism or use of memory or both [4-6].

Any successful metaheuristic algorithm requires a good balance of these two important, seemingly opposite, components [6]. If the intensification is too strong, only a fraction of local space might be visited, and there is a risk of being trapped in a local optimum, as it is often the case for the gradient-based search such as the classic Newton-Raphson method. If the diversification is too strong, the algorithm will converge too slowly with solutions jumping around some potentially optimal solutions. Typically, the solutions start with some randomly generated, or educated guess, solutions, and gradually reduce their diversification while increase their intensification at the same time, though how quick to do so is an important issue.

Another important feature of modern metaheuristics is that an algorithm is either trajectory-based or population-based. For example, simulated annealing is a good example of trajectory-based algorithm because the path of the active search point (or agent) forms a Brownian motion-like trajectory with its movement towards some attractors. On the other hand, genetic algorithms are a good example of population-based method since the parameter search is carried out by multiple genes or agents in parallel. It is difficult to decide which type of method is more efficient as both types work almost equally successfully under appropriate conditions. There are some hints from the recent studies that population-based algorithms might be more efficient for multiobjective multimodal optimization problems as multiple search actions are in parallel, this might be true from the implementation point of view; however, there is far from conclusive and there is virtually no theoretical research to back this up. It seems again that a good combination of these two would lead to better metaheuristic algorithms.

We will now brief introduce some other popular metaheuristic algorithms, and we will try to see how intensification and diversification are used.

## *3.2 Popular Metaheuristic Algorithms*

### 3.2.1 Simulated Annealing

Simulated annealing is probably the best example of modern metaheuristic algorithms, and it was developed by Kirkpatrick, Gelatt and Vecchi in 1983 [7], inspired by the annealing process of metals during heat treatment and Metropolis algorithms for Monte Carlo simulations. The basic idea of the simulated annealing algorithm is similar to dropping some bouncing balls over a landscape, and as the balls bounce and loose energy, they will settle down at some local minima. If the balls are allowed to bounce enough times and loose energy slowly enough, some of the balls will eventually fall into the globally lowest locations, and hence the minimum will be reached. Of course, we can use many balls (parallel simulated annealing), or use a single ball to trace its trajectory (standard simulated annealing).

The optimization process typically starts from an initial guess with higher energy. It then moves to other locations randomly with slightly reduced energy. The move is accepted if the new state has lower energy and the solution improves with a better objective or lower value of the objective function for minimization. However, if the solution does not improve, it is still accepted with a probability of

$$p = \exp[-\frac{\delta E}{kT}], \qquad (8)$$

which is a Boltzmann-type probability distribution. Here $T$ is the temperature of the system, while $k$ is the Boltzmann constant and can be taken to be 1 for simplicity. The energy difference $\delta E$ is often related to the objective function $f(x)$ to be optimized. The trajectory in the simulated annealing is a piecewise path, and this is virtually a Markov chain as the new state (new solution) only depends on the current state/solution and the transition probability $p$.

Here the diversification via randomization produces new solutions (locations), whether the new solution is accepted or not is determined by a probability. If $T$ is too low ($T \to 0$), then any $\delta E > 0$ (worse solution) will rarely be accepted as $p \to 0$, and the diversity of the solutions is subsequently limited. On the other hand, if $T$ is too high, the system is at a high-energy state, most new changes will be accepted, and the minima are not easily reached. So the temperature $T$ is essentially controlling the balance of diversification and intensification. The change of $T$ is called the cooling schedule.

There are two main categories of cooling schedules: the monotonically decrease and non-monotonic. For monotonic cooling, geometric cooling schedule is by far the most widely used: $T(t) = T_0 \alpha^t$, where $t$ is the time step or counter of iterations and $T_0$ is the initial temperature. The advantage of this schedule is that there is no need to determine the final temperature, and $\alpha$ is in the range of (0,1). The disadvantage is that if you use a very small value of $\alpha$ corresponding to the simulated quenching (SQ), then there is a risk for the system to freeze too quickly and it might be trapped in some local optima. If $\alpha$ is approaching 1, then the convergence is slow. A possible solution is to use non-monotonic cooling schedule so that the system can be elevated to a higher energy state when necessary. Again, if you raise the temperature too many times, the convergence is affected. This demonstrates that there is a fine balance between diversification and intensification.

### 3.2.2 Evolutionary Algorithms

Evolutionary algorithms (EA) are the name for a subset of evolutionary compution [8] and they are the search methods which are inspired from Charles Darwin's natural selection and survival of the fittest. Evolutionary algorithms are population-based search algorithms and almost all of them use genetic operators to a certain degree. These operators typically include crossover or reproductive recombination, mutation, inheritance and selection based on their fitness.

Genetic algorithms (GA) are by far the most popular and widely used [9], and other evolutionary search methods are almost equally popular, and these include genetic programming, evolutionary programming, and evolutionary strategies. Briefly speaking, genetic programming (GP) is an evolutionary machine-learning technique in the framework of genetic algorithms where each individual in the population is a computer program using a scheme-style computer language such as Lisp. The objective is to find the optimal computer program to perform a user-defined task.

Evolutionary strategy (ES) is another class of nature-inspired evolutionary optimization techniques, and it mainly uses mutation, selection and element-wise average for intermediate recombination as the genetic operators. It has been applied to a wide range of optimization problems.

Evolutionary programming (EP) uses arbitrary data structures and representations tailored to suit a particular problem domain, and they are combined with the essence of genetic algorithms so as to solve generalized complex optimization problems. EP is very similar to ES, but without the recombination operator, and the main difference between EP and other methods is that EP does not use exchange of string segments, and thus there is no crossover or recombination between individuals. The primary genetic operator is mutation, often using Gaussian mutation for real-valued functions. However, its modern variants are more diverse. Here we will briefly introduce the basic idea of genetic algorithms.

Genetic algorithms were developed by John Holland in the 1960s and 1970s [9], using crossover and recombination, mutation and selection for adaptive and artificial systems, optimization and other problem-solving strategies. The essence of genetic algorithms involves the encoding of an optimization function or a set of functions into arrays of binary or character strings to represent chromosomes, and the manipulation of these strings by genetic operators with the ultimate aim to find the optimal solution for the problem concerned. The encoding is often carried out into multiple binary strings, called population, though real-valued encoding and representations are more often used in modern applications. The initial population then evolves by generating new-generation individuals via crossover of two randomly selected parent strings, and/or the mutation of some random bits. Whether a new individual is selected or not is based on its fitness, which is linked in some way to the objective functions.

The frequency of crossover represented by a crossover probability $p_c$ in the range [0,1]. Similarly, the mutation is controlled by a mutation rate $p_m$. The crossover of two parent strings is the main operator with a higher probability, often in the range of 0.7 to 0.99, and is carried out by swapping one segment of one chromosome with their corresponding segment on another chromosome at a random position or even multiple positions. If the crossover probability is too low, crossover occurs sparsely, which means the population changes slowly, and thus not efficient to explore all possibilities and the part of the search space it explores is also limited.

Mutation is often carried out by flipping some bits on a chromosome, and the mutation probability $p_m$ is typically low in the range of 0.001 to 0.1. If the mutation rate is too high, the individuals produced might be too different from their main population, and the solutions could 'jump around' some optima without settling down or converging. Ideally, when the optimal solution is approaching, the mutation rate should be reduced gradually so that a good rate of convergence can be achieved. Thus, a good balance of the rate of mutation and the rate of crossover should be maintained.

Another important issue is the selection of individuals based on their fitness. There are many ways to define the fitness function, which has to be linked to the objective function f($x$). A popular and yet simple example is to use the relative fitness $F_i$ for the whole population

$$F_i = \frac{f(\xi_i)}{\sum_{i=1}^{N} f_i(\xi_i)} \qquad (9)$$

where $\xi_i$ is the phenotypic value of the individual *i*, and it corresponds to the solution vector in most applications. Here *N* is the population size, and the right choice is also important. If the population size is too small, there is not enough evolution going on, and there is a risk for the whole population to go extinct in nature. Numerically, this means that whole population set may go astray or be dominated by a few individuals, leading to premature convergence at meaningless solutions. On the other hand, if the population is too large, the diversity of the solution might be too extreme, resulting in slow convergence and more computing time.

Another subtle issue is the use of elitism in the selection of the fittest. The simplest elitism is to select a small set of the most fit individuals in each generation which will be carried over to the new generation without being modified by any genetic operators. This essentially ensures that the best solutions always remain in the population, and subsequently a quicker convergence with a guaranteed quality set of best solutions may be achieved. There exists an associated issue here and that is the extent of the elitism. If the elitism is too extensive, most individuals even the not-so-good solutions are carried over, the search space is not well explored and the quality of the solutions may be affected and the convergence is slow if it ever converges. On the other hand, if there is little elitism, it is not much different from the standard genetic algorithms, and the use of elitism might have no effect in improving the convergence at all. Again some fine balance might be needed.

The advantages of the genetic algorithms over traditional optimization algorithms are the ability of dealing with complex problems and parallelism. GA search methods can deal with various types of optimization, whether the objective functions are stationary or non-stationary (time-dependent), linear or nonlinear, continuous or discontinuous, or even with random noise. As individuals in a population act like independent agents, each subset of the population can explore the search space in many directions simultaneously. This feature makes it ideal for the parallel implementation of the genetic algorithms.

### 3.2.3 Particle Swarm Optimization

As simulated annealing is a trajectory-based metaheuristic and genetic algorithms are population-based, we now introduce another population-based metaheuristic algorithm, namely, particle swarm optimization (PSO). The PSO was developed by Kennedy and Eberhart in 1995 [10], inspired by the swarm behaviour of fish and bird schooling in nature. Unlike the single trajectory used in simulated annealing, this algorithm searches the space of the objective function by adjusting the multiple trajectories of individual agents, called particles. The motion of the particles has two major components: a stochastic component and a deterministic component in terms of velocity and position vectors (solution vectors)

$$v_i^{t+1} = v_i^t + \alpha\, \varepsilon_1 (x_i - g^*) + \beta\, \varepsilon_2 (x_i - x_i^*), \qquad x_i^{t+1} = x_i^t + v_i^t, \qquad (10)$$

where $v_i^t$ and $x_i^t$ are the velocity and position of particle *i* at time *t*, respectively. $\varepsilon_1$ and $\varepsilon_2$ are two random vectors, while $\alpha$ and $\beta$ are constants, often called the learn-ing parameters. The diversification is controlled by the combination of random vectors and the learning or diversification parameters.

The intensification is mainly represented by the deterministic motion towards the updated current best $x_i^t$ for particle $i$, and the current global best $g^*$ for all particles. As the particles approach to the optima, their motion and randomness are reduced. There are many different variants of PSO in the literature [4,10].

There is a hidden or implicit feature in the PSO algorithm, that is the broadcasting ability of the current global best $g^*$ to other particles. That can be thought as either the use of memory or some higher-level strategy so as to speed up the convergence and explore the search space more effectively and efficiently. If the diversification parameter is large, a larger part of the search space will be explored; however, it will converge more slowly. On the other hand, high-level intensification will make the algorithm converge quickly, but not necessarily to the right solution set. Also if the intensification is too strong, it may even slow down the convergence rate as the system is still slowly evolving. Again, a fine balance between the diversification and intensification is important.

### 3.2.4 Ant Colony Optimization

Another population-based metaheuristic algorithm is the ant colony optimization (ACO) which was first formulated by Dorigo and further developed by other pioneers [11-13]. This algorithm was based on the characteristics of behaviour of social ants. For discrete routing and scheduling problems, multiple ants are often used. Each virtual ant will preferably choose a route with higher pheromone concentration, and it also deposits more pheromone at the same time. If there is no previously deposited pheromone, then each ant will move randomly. In addition, the pheromone concentration will decrease gradually due to the evaporation, often with a constant evaporation rate.

Here the diversification of the solutions is represented by the randomness and the choice probability of agents along a specific route. The intensification is implicitly manipulated by the pheromone concentration and the evaporation rate. However, the evaporation rate can also affect the diversification in some way. This algorithm is exceptionally successful in the combinatorial optimization problems. Again, some fine balance between the diversification and intensification is needed to ensure a faster and efficient convergence and to ensure the quality of the solutions.

### 3.2.5 Firefly Algorithm

The fascinating flashing light of fireflies in the tropical summer can be used to develop interesting nature-inspired metaheuristic algorithms for optimization. The Firefly Algorithm (FA) was developed by Xin-She Yang in 2007 [4], based on the idealization of the flashing characteristics of fireflies. There are three major components in the FA optimization: 1) A firefly will attract to more brighter or more attractive fireflies, and at the same time they will move randomly; 2) the attractiveness is proportional to the brightness of the flashing light which will decrease with distance, therefore, the attractiveness will be evaluated in the eye of the beholders (other fireflies); 3) The decrease of light intensity is controlled by the light absorption coefficient γ which is in turn linked to a characteristic scale.

The new solution is generated by

$$x_i^{t+1} = x_i^t + \alpha\, \varepsilon_1 + \beta \exp[-\gamma\, r_{ij}^2]\, \varepsilon_2 (x_i - x_j), \qquad (10)$$

where $r_{ij}$ is the distance, not necessarily the physical distance, between two fireflies $i$ and $j$. In addition, FA is obviously a population-based algorithm, which may share many similarities with particle swarm optimization. In fact, it has been proved by Yang [4] that when $\gamma \to \infty$, the firefly algorithm will become an accelerated version of PSO, while $\gamma \to 0$, the FA reduces to a version of random search algorithms.

In the FA optimization, the diversification is represented by the random movement component, while the intensification is implicitly controlled by the attraction of different fireflies and the attractiveness strength β. Unlike the other metaheuristics, the interactions between exploration and exploitation intermingled in some way; this might be an important factor for its success in solving multiobjective and multimodal optimization problems. As we will discuss later, the hybridization of diversification and intensification components is a useful technique to development new algorithms. Again, a good balance between these two factors is very important.

Obviously, there are many other metaheuristic algorithms which are currently widely used, including, but not limited to, evolutionary algorithms, artificial bee algorithms [14], Tabu search, photosynthetic algorithm, enzyme algorithm [15], cross-entropy algorithm, and last but not least, the harmony search we discussed earlier. Readers can find more references in [4-9,13-19]. After the brief introduction to other metaheuristic algorithms, we are now ready to analyze the similarities and differences of the Harmony Search algorithm in the general context of metaheuristics.

## 4 Characteristics of HS and Comparison

### *4.1 Diversification and Intensification*

In reviewing other metaheuristic algorithms, we have repetitively focused on two major components: diversification and intensification. They are also referred to as exploration and exploitation [6,13]. These two components are seemingly contra-dicting each other, but their balanced combination is crucially important to the success of any metaheuristic algorithms [4,6,13].

Proper diversification or exploration makes sure the search in the parameter space can explore as many locations and regions as possible in an efficient and effective manner. It also ensures that the evolving system will not be trapped in biased local optima. Diversification is often represented in the implementation as the randomization and/or additional stochastic component superposed onto the deterministic components. If the diversification is too strong, it may explore more accessible search space in a stochastic manner, and subsequently will slow down the convergence of the algorithm. If the diversification is too weak, there is a risk that the parameter space explored is so limited and the solutions are biased and trapped in local optima, or even lead to meaningless solutions.

On the other hand, the appropriate intensification or exploitation intends to exploit the history and experience of the search process. It aims to ensure to speed up the convergence when necessary by reducing the randomness and limiting diversification. Intensification is often carried out by using memory such as in Tabu search and/or elitism such as in the genetic algorithms. In other algorithms, it is much more elaborate to use intensification such as the case in simulated annealing and firefly algorithms. If the intensification is too strong, it could result in premature convergence, leading to biased local optima or even meaningless solutions, as the search space is not well explored. If the intensification is too weak, convergence becomes slow.

The optimal balance of diversification and intensification is required, and such a balance itself is an optimization process. Fine-tuning of parameters is often required to improve the efficiency of the algorithms for a particular problem. There is No Free Lunch in any optimization problems [18]. A substantial amount of studies might be to choose the right algorithms for the right optimization problems [19], though it lacks a systematic guidance for such choices.

## *4.2 Why HS is Successful*

Now if we analyze the Harmony Search algorithm in the context of the major components of metaheuristics and try to compare with other metaheuristic algorithms, we can identify its ways of handling intensification and diversification in the HS method, and probably understand why it is a very successful metaheuristic algorithm.

In the HS algorithm, diversification is essentially controlled by the pitch adjustment and randomization -- here there are two subcomponents for diversification, which might be an important factor for the high efficiency of the HS method. The first subcomponent of composing 'new music', or generating new solutions, via randomization would be at least at the same level of efficiency as other algorithms by randomization. However, an additional subcomponent for HS diversifi-cation is the pitch adjustment characterized by $r_{pa}$. Pitch adjusting is carried out by adjusting the pitch in the given bandwidth by a small random amount relative to the existing pitch or solution from the harmony memory. Essentially, pitch adjusting is a refinement process of local solutions. Both memory considering and pitch adjusting ensure that the good local solutions are retained while the randomization and harmony memory considering will explore the global search space effectively. The subtlety of this is that it is a controlled diversification around the good solutions (good harmonics and pitches), and it almost acts like an intensification factor as well. The randomization explores the search space more efficiently and effectively; while the pitch adjustment ensures that the newly generated solutions are good enough, or not too far way from existing good solutions.

The intensification is mainly represented in the HS algorithm by the harmony memory accepting rate $r_{accept}$. A high harmony acceptance rate means the good solutions from the history/memory are more likely to be selected or inherited. This is equivalent to a certain degree of elitism. Obviously, if the acceptance rate is too low, the solutions will converge more slowly. As mentioned earlier, this intensification is enhanced by the controlled pitch adjustment. Such interactions between various components could be another important factor for the success of the HS algorithm over other algorithms, as it will be demonstrated again and again in later chapters in this book.

In addition, the implementation of HS algorithm is also easier. There is some evidence to suggest that HS is less sensitive to the chosen parameters, which means that we do not have to fine-tune these parameters to get quality solutions.

Furthermore, the HS algorithm is a population-based metaheuristic, this means that multiple harmonics groups can be used in parallel. Proper parallelism usually leads to better implantation with higher efficiency. The good combination of parallelism with elitism as well as a fine balance of intensification and diversification is the key to the success of the HS algorithm, and in fact, to the success of any metaheuristic algorithms.

These advantages make it very versatile to combine HS with other metaheuristic algorithms such as PSO to produce hybrid metaheuristics [20] and to apply HS in various applications [1-3,20-22].

## 5 Further Research

The power and efficiency of the HS algorithm seem obvious after the discussion and comparison with other metaheuristics; however, there are some unanswered questions concerning the whole class of HS algorithms. At the moment, the HS algorithm like almost all other metaheuristics is a higher-level optimization strategy which works well under appropriate conditions, but we usually do not fully understand why and how they work so well. For example, when choosing the harmony accepting rate, we usually use a higher value, say, 0.7 to 0.95. This is obtained by experimenting the simulations, or using a similar inspiration from ge-netic algorithms when the mutation rate should be low, and thus the accepting rate of the existing gene components are high. However, it is very difficult to say what range of values and which combinations are surely better than others.

In general, there lacks a theoretical framework for metaheuristics to provide some analytical guidance to the following important issues: How to improve the efficiency for a given problem? What conditions are needed to guarantee a good rate of convergence? How to prove the global optima are reached for the given metaheuristic algorithm? These are still open questions that need further research. The encouraging thing is that many researchers are interested in tackling these difficult challenges, and important progress has been made concerning the convergence of algorithms such as simulated annealing. Any progress concerning the convergence of HS and other algorithms would be influentially profound.

Even without a solid framework, this does not discourage scientists to development more variants and/or hybrid algorithms. In fact, the algorithm development itself is a metaheuristic process in a similar manner to the key components of HS algorithms: to use the existing successful algorithms; to develop slightly different variants based on the existing algorithms; and to formulate heuristically completely new metaheuristic algorithms. By using the existing algorithms, we have to found the right algorithms for the right problems. Often, we have to change and reformulate the problem slightly and/or to improve the algorithms slightly so as to find the solutions more efficiently. Sometimes, we have to develop new algorithms from scratch to solve some tough optimization problems.

There are many ways to develop new algorithms, and from the metaheuristic point of view, the most heuristic way is probably to develop new algorithms by hybridization. That is to say, new algorithms are often based on the right combination of the existing metaheuristic algorithms. For example, combining a trajectory-type simulated annealing with multiple agents, the parallel simulated annealing can be developed. In the context of HS algorithms, the combination of HS with PSO, the global-best harmony search has been developed [20]. As in the case of any efficient metaheuristic algorithms, the most difficult thing is probably to find the right or optimal balance between diversity and intensity of the found solutions; here the most challenging task in developing new hybrid algorithms is probably to find the right combination of which feature/components of existing algorithms.

A further extension of the HS algorithm will be to solve multiobjective optimization problems more naturally and more efficiently. At the moment, most of the existing studies, though very complex and tough *per se*, have been mainly focused on the optimization with a single objective or a few criteria with linear and nonlinear constraints. The next challenges would be to use the HS algorithm to solve tough multiobjective and multicriterion NP-hard optimization problems.

Whatever the challenges will be, more and more metaheuristic algorithms will be applied in various applications [21-25] as well as more systematic studies [26-27]. Subsequently more and more new hybrid metaheuristic algorithms will be developed in the future. Furthermore, more solid theoretical work will pace the way for further research and provide some guidance to new algorithm formulations. As you will see the rest chapters of the book, metaheuristic algorithms such as Har-mony Search will show again and again their power, novelty, effectiveness and efficiency in solving many tough optimization problems in a diverse range of applications in sciences, engineering and industry.

## References

1. Geem ZW, Kim JH and Loganathan GV (2001) A new heuristic optimization algorithm: Harmony search. Simulation, 76:60-68
2. Lee KS and Geem ZW (2005) A new meta-heuristic algorithm for continuous engineering optimization: harmony search theory and practice. Comput. Methods Appl. Mech. Engrg., 194:3902-3933
3. Geem ZW (2006) Optimal cost design of water distribution networks using harmony search. Engineering Optimization 38:259-280
4. Yang XS (2008) Nature-inspired Metaheuristic Algorithms. Luniver Press
5. Glover F and Laguna M (1997) Tabu Search. Kluwer Academic Publishers
6. Blum C and Roli A (2003) Metaheuristics in combinatorial optimization: Overview and conceptual comparison. ACM Comput. Surv., 35:268-308
7. Kirkpatrick S, Gelatt CD and Vecchi MP (1983) Optimization by simulated annealing. Science, 220:671-680
8. Back T, Fogel D and Michalewicz Z (1997) Handbook of Evolutionary Computation. Oxford University Press
9. Holland JH (1975) Adaptation in Natural and Artificial Systems. The University of Michigan Press, Ann Arbor
10. Kennedy J and Eberhart RC (1995) Particle swarm optimization. Proceedings of IEEE Int. Conf. Neural Networks, pp.1942-1948
11. Dorigo M and Stutzle T (2004) Ant Colony Optimization. MIT Press, Cambridge
12. Bonabeau E, Dorigo M and Theraulaz G (1999) Swarm Intelligence:From Natural to Artificial Systems. Oxford University Press
13. Dorigo M and Blum C (2005) Ant colony optimization theory: A survey. Theor. Comput. Sci., 344:243-278
14. Karaboga D and Basturk B (2008) On the performance of artificial bee colony (ABC) algorithm. Applied Soft Computing, 8:687-697
15. Yang XS (2005) New enzyme algorithm, Tikhonov regularization and inverse parabolic analysis. Advances in Computational Methods in Science and Engineering – ICCMSE'05 Eds T Simos and G Maroulis, 4:1880-1883
16. Yang XS (2005) Biology-derived algorithms in engineering optimization. Handbook of Bioinspired Algorithms and Applications, Eds Olarius S and Zomaya A, Chapman & Hall/CRC
17. Engelbrecht AP (2005) Fundamentals of Computational Swarm Intelligence. Wiley
18. Wolpert DH and Macready WG (1997). No free lunch theorems for optimization. IEEE Transaction on Evolutionary Computation, 1:67-82
19. Yang XS (2008) Mathematical Optimization: From Linear Programming to Metaheuristics, Cambridge Int. Science Publishing, UK.
20. Omran M and Mahdavi (2008) Global-best harmony search. Applied Math. Computation, 198:643-656

21. Geem ZW (2008) Harmony search for solving Sudoku, Knowledge-Based Intelligent Information and Engineering Systems, LNAI, 4692:371-378
22. Geem ZW (2007) Optimal scheduling of multiple dam system using harmony search algorithm. Lecture Notes in Computer Science 4507:316-323
23. Perelman L, Ostfeld A (2007) An adaptive heuristic cross-entropy algorithm for optimal design of water distribution systems. Engineering Optimization 39:413-428
24. Chatterjee A and Siarry P (2006) Nonlinear inertia variation for dynamic adaptation in particle swarm optimization. Comp. Oper. Research, 33:859-871
25. Jaeggi D, Parks GT, Kipouros T, Clarkson PJ (2005) A multiobjective Tabu search algorithm for constrained optimized problems, 3rd Int. Conf. Evolutionary Multi-criterion Optimization, Mexico, LNCS, 3410:490-504
26. Michalewicz Z (1996) Genetic Algorithm + Data Structure = Evolutionary Programming, Springer, New York
27. Spall JC (2003) Introduction to Stochastic Search and Optimization: Estimation, Simulation and Control, Wiley.